\documentclass[a4paper, 12pt]{amsart}

\usepackage [latin1]{inputenc}
\usepackage[all]{xy}
\usepackage{amsmath}
\usepackage{amsfonts}
\usepackage{latexsym}  
\usepackage{amssymb}
\usepackage{graphics}
\usepackage{graphicx}   
\newtheorem{teorema}{Theorem}
\newtheorem{lemma}[teorema]{Lemma}
\newtheorem{propos}[teorema]{Proposition}
\newtheorem{corol}[teorema]{Corollary}
\newtheorem{ex}[teorema]{Example}
\newtheorem{rem}[teorema]{Remark}
\newtheorem{defin}[teorema]{Definition}

\def\defin{\par\ifdim\lastskip<\smallskipamount\removelastskip
  \smallskip\fi\noindent{\bf\ignorespaces
Definition\unskip:\enspace}\rm \ignorespaces}

\def\bit{\begin{itemize}}
\def\eit{\end{itemize}}
\def\be{\begin{equation}}
\def\ee{\end{equation}}
\def\beq{\begin{eqnarray}}
\def\eeq{\end{eqnarray}}

\def\ba{\begin{array}}
\def\ea{\end{array}}
\def\bt{\begin{teorema}}
\def\et{\end{teorema}}
\def\bp{\begin{propos}}
\def\ep{\end{propos}}
\def\bl{\begin{lemma}}
\def\el{\end{lemma}}
\def\bc{\begin{corol}}
\def\ec{\end{corol}}
\def\br{\begin{rem}\rm}
\def\er{\end{rem}}
\def\bex{\begin{ex}\rm}
\def\eex{\end{ex}}
\def\bd{\begin{defin}}
\def\ed{\end{defin}}
\def\demo{\par\noindent{\bf Proof.\ }}
\def\enddemo{\vspace {-0.55cm}
\begin{flushright}
$\Box$
\end{flushright}
\par\vskip.6truecm}

\def\ES{\emptyset}    

\def\R{{\mathbb {R}}}   \def\a {\alpha} \def\b {\beta}
\def\N{{\mathbb N}}     \def\d {\delta} \def\e{\varepsilon}
\def\C{{\mathbb C}}      
					\def\o{\omega}\def\p{\partial}
					\def\t{\theta}\def\z{\zeta}\def\n{\nu}
\def\P{{\mathbb P}}					

																								\def\G{\Gamma}
																								\def\O{\Omega}					

  \def\smi{\smallsetminus} \def\ssmi{\!\smallsetminus\!}

\def\sbs{\!\subset\!}

\def\oli{\overline}

\def\benu{\begin{enumerate}} \def\eenu{\end{enumerate}}
\def\beqn{\begin{eqnarray*}}  \def\eeqn{\end{eqnarray*}}

\def\beqn{\begin{eqnarray*}}  \def\eeqn{\end{eqnarray*}}

\def\bd{\rm{b}\,}
\def\dist{\rm{dist}\,}

\begin{document}
\title[The Bremermann-Dirichlet problem for unbounded domains]{The Bremermann-Dirichlet problem for unbounded domains of $\C^n$}
\author{Alexandru Simioniuc, Giuseppe Tomassini}
\address{Scuola Normale Superiore, Piazza dei Cavalieri, 56126-Pisa, Italy}
\email{a.simioniuc@sns.it, g.tomassini@sns.it}

\begin{abstract}

Given an unbounded strongly pseudoconvex domain $\Omega$ and a continuous real valued function $h$ defined on $\bd\O$, we study the existence of a (maximal) plurisubharmonic function $\Phi$ on $\O$ such that $\Phi|_{\bd\O}=h$.


\end{abstract}

\maketitle

\tableofcontents
\section{Introduction}
Let $\O\sbs\C^n$ be a bounded domain and $h:{\rm b}\O\to \R$ a continuous function. The problem of extending $h$ to $\O$ by a plurisubharmonic function was considered for the first time by Bremermann in \cite{Br59}. He proved that, if $\O$ is strongly pseudoconvex, the upper envelope $u_{\O,h}$ of the class of plurisubharmonic functions in $\O$ which are majored by $h$ on ${\rm b}\O$, is a plurisubharmonic extension of $h$ which is continuous at ${\rm b}\O$. Later on, Walsh (\cite{Wa69}) showed that $u_{\O,h}$ is actually continuous and Bedford and Taylor (\cite{BT76}) proved that $u_{\O,h}$ ias a solution of the homogeneous complex Monge-Amp\`ere equation.

It is worth observing that, if $\O$ is not strongly pseudoconvex, the boundary value $h$ cannot be arbitrary (cfr.\ Lemma \ref{l1}). 

A generalization of the problem for $q$-plurisubharmonic functions and strongly $q$-pseudoconvex domains was considered by Hunt and Murray (\cite{Hu78}) and by Slodkowski (\cite{Sl84}).

In this paper we deal with the following {\it Bremermann-Dirichlet problem} consisting of finding a function $u:\oli\O\to\R$, such that\vspace{5mm}
\begin{itemize}
	\item [($\star$)] $\left\{\begin{array}{cl}(i)& u \text{ is upper semicontinuous in }\oli\O\\ (ii)&u\text{ is plurisubharmonic in }\O\\ (iii)&u|_{\bd\O}=h\end{array}\right.$
\end{itemize}
\vspace{5mm}
where $\O\sbs\C^n$ is an unbounded strongly pseudoconvex domain with a $C^2$ boundary and $h:{\rm b} \O\to\R$ a continuous function. 

In this situation the problem may not admit non-trivial solutions even for very simple domains. Indeed, in \cite{ST99}, the authors constructed a continuous function $h$ on the boundary ${\rm b} \O$ of a paraboloid $\O$, satisfying $\inf\limits_{{\rm b} \O}\,h=-\infty$, such that the only plurisubharmonic function $u$ satisfying $u\le h$ on ${\rm b} \O$ is the function $u\equiv-\infty$. 

Thus, in ($\star$) we assume that $h\geq 0$.

In the first part of the paper we prove that the problem ($\star$) admits a solution $u$ which is continuous in $\oli\O$ (cfr.\ Theorem \ref{constr}).

In the second part of the paper we study the existence of a maximal solution for the Bremermann-Dirichlet problem. We first consider the case of strongly convex domains. Using an appropriate exhaustive sequence of subdomains of $\O$ and a decreasing sequence of Perron-Bremermann functions defined on bounded domains we prove the following:
\begin{itemize}
\item Let $\O\subset\C^n$ be an unbounded strongly convex domain and $h:\bd\O\to\R$ a bounded continuous function. Then the Bremermann-Dirichlet problem admits a maximal solution, which is continuous on $\oli\O$.
\end{itemize}
(cfr.\ Proposition \ref{2345}). 
If $h:\bd\O\to[0,+\infty)$ is not bounded, the continuity of the maximal solution, which is always granted on the boundary $\bd \O$ (cfr.\ Proposition \ref{6789}) is hard to prove in general. In Subsection \ref{se} we exhibit some examples where we can obtain the continuity of the maximal solution assuming some hypothesis on the boundary of the domain or on the function $h$.

The general case of a strongly pseudoconvex domain is treated in Section \ref{9999}. The existence of the maximal solution and its continuity at the boundary are obtained assuming an aditional hypothesis for $\O$, which was introduced by Lupacciolu (\cite{Lu87}) in studying the extension problem for $CR$ functions in unbounded domains.

The last section of the paper deals with $q$-plusubharmonic solutions. We prove that if $\O$ is an unbounded strongly convex domain and $h$ a bounded real valued continuous function defined on $\bd\O$, then the generalized Bremermann-Dirichlet problem for $q$-plusubharmonic functions admits a maximal solution which is continuous on $\oli\O$ and ($n-q-1$)-plurisuperharmonic in $\O$.
\section{Preliminaries}
\subsection{The Perron-Bremermann function}
A plurisubharmonic function $u:U\to\R\cup\{-\infty\}$ on an open set $U\subseteq\C^n$ is an upper semicontinuous function such that its restriction to any complex line $L$ is subharmonic, i.e.\ for every compact subset $K\subset U\cap L$ and for every harmonic function $h$ on $\rm{int}\,K$ continuous on $K$ such that $h|_{\bd K}\geq u|_{\bd K}$ we have $h\geq u|_{K}$. If $U\in\C^n$ is an open set, we denote by ${\sf Psh}(U)$ the set of plurisubharmonic functions in $U$ and by ${\sf Psh}^{(s)}(U)$ the subset of functions $u\in{\sf Psh}(U)$ which are semicontinuous on $\oli U$.
\bp\label{lmax}
Let $U, V\subset\C^n$ be open sets, $V\subset U$, $u\in {\sf Psh}^{(s)}(U)$ and $v\in {\sf Psh}^{(s)}(V)$ such that $u|_{\bd V}\geq v|_{\bd V}$. The function $\Psi:\oli U\to\R\cup\{-\infty\}$ given by
$$
\Psi(z)=\left\{\begin{array}{ll} \max\{u(z),v(z)\},& z\in\oli V\\ u(z),& z\in\oli U\ssmi\oli V,\end{array}\right.
$$
is plurisubharmonic on $U$.
\ep
\demo
Let $m=\max\{u,v\}$. It is clear that $\Psi$ is upper semicontinuous. Let $L$ be a complex line, $K\subset L\cap\oli U$ a compact subset and $h:K\to\R$ a continuous function, harmonic on $\rm{int}\,K$, such that $h|_{\bd K}\geq\Psi|_{\bd K}$. We have to prove that $h\geq\Psi|_{K}$. Suppose by contradiction that there exists $z^0\in K$ such that $h(z^0)<\Psi(z^0)$. Since $\Psi$ is plurisubharmonic on $V$ and on $U\ssmi \oli V$, we must have $K\cap\bd V\neq\ES$. If $z^0\in V$ (respectively $z^0\in U\ssmi\oli V$), by the maximum principle applied to the subharmonic function $m|_L$ on $V\cap L$ (repectively $u|_L$ on $U\cap L\ssmi\oli V$), we deduce in both cases that there exists $z^1\in K\cap \bd V$ such that $$u(z^1)=m(z^1)>h(z^1).$$
This contradicts the maximum principle for $u|_L$, since $h|_{\bd K}\geq\Psi|_{\bd K}\geq u|_{\bd K}$.
\enddemo

A plurisubharmonic function $u:\O\to\R\cup\{-\infty\}$ is called {\it(locally) maximal} if for every open set $G\subset\O$ and every function $v$ plurisubharmonic on $G$ such that $\limsup_{z\to p} v(z)\leq u(p)$ for all $p\in\bd\O$ we have $u\leq v$ on $G$.

Let $\O\subseteq\C^n$ be a domain with $C^2$ boundary. Let $p\in\bd \O$ and $U$ be an open neighbourhood of $p$ in $\C^n$. A $C^2$ function $\varrho:U\to\R$ satisfying
$$D\cap\O=\{z\in U:\ \varrho(z)<0\}$$
and $d\varrho\neq0$ on $\bd\O\cap U$ is called a {\it defining function} on $U$ for ${\rm b}\O$. 

Let $L_p(\varrho)$ denote the Levi form 
$$\oli\p\p \varrho(p)=\sum_{j,k=1}^n\frac{\p^2\varrho(p)}{\p z_j\p\oli z_k}dz_i\wedge d\oli z_j$$
of $\varrho$ at $p$. $\O$ is called strongly pseudoconvex at $p$ if the restriction of $L_p(\varrho)$ to the tangent hyperplane to $\bd\O$ at $p$ is positively definite. The definition does not depend on the defining function $\varrho$. A bounded strongly pseudoconvex domain $\O$ admits a plurisubharmonic globaly defining function in a neighbourhood of $\oli\O$, i.e.\ can be written in the form
$$\O=\{z\in\C^n:\varrho(z)<0\},$$
where $u$ is strongly plurisubharmonic in a neighborhood of $\oli\O$ and ${\rm d}\varrho\neq 0$ on $\bd\O$.

We state here, the following
\bl \label{l1}
Let $\O\in\C^n$, $n\geq 2$, be a domain and $H$ a hyperplane such that $H\cap\O\neq\emptyset$. Let $\O_1$ be a bounded connected component of $\O\smi H$ and $\Sigma$ and $L$ such that $\bd\O_1=\Sigma\cup L$, where $L\neq\emptyset$ is an open subset of $H$. If $u\in{\sf Psh}^{(s)}(\O)$, then
$$\max_{z\in\oli\O_1}u(z)=\max_{z\in \Sigma}u(z).$$
\el
\demo
Let $z^0\in\oli \O_1$. We choose a complex line $L_1\subset \oli\O_1$ passing by $z^0$ such that $L_1\cap \O$ is a domain $D$ of $L_1$ with $\bd D\subset \Sigma$. Considering the constant function $\phi\equiv M$ on $\oli D$, where $M=\max_{z\in \Sigma}u(z)$, $\phi$ is a harmonic function satisfying $u|_{bd D}\leq \phi|_{bd D}$, therefore $u|_{D}\leq\phi$ and thus $u(z^0)\leq M$.
\enddemo
For every continuous function $h:\bd\O\to\R$, we define
$$
\mathcal P_{\O,h}=\{u\in {\sf Psh}^{(s)}(\O):\ u|_{\bd\O}\leq h\}
$$
and we denote $u_{\O,h}:\oli\O\to [-\infty,+\infty)$ the upper envelope of $\mathcal P_{\O,h}$
$$
u_{\O,h}(z)=\limsup_{\z\to z}\sup_{u\in\mathcal P_{\O,h}}u(\z).
$$
$u_{\O,h}$ is called the {\it Perron-Bremermann function} for $\O$ and $h$.
\bp \label{p1}
Let $\O$ be a bounded domain in $\C^n$. The Perron-Bremermann function has the following properties:
\begin{itemize}
	\item [(i)] $u_{\O,h+c}=u_{\O,h}+c,\ \forall c\in\R$;
	\item [(ii)] if $\ h_1\leq h_2$, then $\ u_{\O,h_1}\leq u_{\O,h_2}$;
	\item [(iii)] if $\ h=\sup(h_1,h_2)$, then $\ u_{\O,h}\geq \sup(u_{\O,h_1},u_{\O,h_2})\ $ (and in general $u_{\O,h}\neq \sup(u_{\O,h_1},u_{\O,h_2}))$;
	\item [(iv)] $\max_{z\in\oli\O}|u_{\O,h_1}-u_{\O,h_2}|=\max_{z\in \bd\O}|h_1-h_2|;$
	\item[(v)] if $h_\nu\to h$ in $C^0(\bd\O)$, then $u_{\O,h_\nu}\to u_{\O,h}$ in $C^0(\oli\O).$
\end{itemize}
\ep
\demo
(i), (ii) and (iii) are obvious. To prove (iv), we put
$$M=\max_{z\in \bd\O}(h_1(z)-h_2(z))$$
and
$$m=\min_{z\in \bd\O}(h_1(z)-h_2(z)).$$
Using (i), (ii) and (iii) we get:
$$-\max(|m|,|M|)\leq m\leq |u_{\O,h_1}-u_{\O,h_2}|\leq M\leq\max(|m|,|M|)$$
and from this it follows the conclusion, since
$$\max(|m|,|M|)=\max_{z\in \bd\O}|h_1-h_2|.$$
Property (v) follows immediately from (iv).
\enddemo
\subsection{The Bremermann-Dirichlet problem}
Given a domain $\O\subset\C^n$ and a continuous function $h:\bd\O\to\R$, the {\it Bremermann-Dirichlet problem} consists of finding a function $u:\oli\O\to\R$, such that\vspace{5mm}
\begin{itemize}
	\item [($\star$)] $\left\{\begin{array}{cl}(i)& u \text{ is upper semicontinuous in }\oli\O\\ (ii)&u\text{ is plurisubharmonic in }\O\\ (iii)&u|_{\bd\O}=h\end{array}\right.$
\end{itemize}
\vspace{5mm}
Bremermann proved in \cite{Br59} that if $\O$ is a strongly pseudoconvex bounded domain, then the Perron-Bremermann function $u_{\O,h}$ is an element of $\mathcal P_{\O,h}$, which is continuous on $\bd \O$ and $u_{\O,h}|_{\bd\O}=h$. In \cite{Wa69}, Walsh proved that $u_{\O,h}$ is actually continuous in $\oli\O$. Thus, for a strongly pseudoconvex bounded domain $\O$, $u_{\O,h}$ is a continuous maximal solution of the problem ($\star$).  

About regularity we have the following
\bp\label{lip1}
Let $\O\subset\C^n$ be a strongly pseudoconvex bounded domain and $h:\bd\O\to\R$ a $C^2$ function. Then the Perron-Bremermann solution $u_{\O,h}$ of the problem ($\star$) is Lipschitz in $\oli\O$.
\ep
We first prove the 
\bl\label{lip}
Let $\O\subset\C^n$ be a strongly pseudoconvex bounded domain and $h:\bd\O\to\R$ a continuous function. If the Perron-Bremermann solution $u_{\O,h}$ of the problem ($\star$) is Lipschitz in a $\d$-neighbourhood $\O_\d$ of $\bd\O$ in $\oli\O$, for some $\d>0$, then $u_{\O,h}$ is Lipschitz in $\oli\O$. Moreover, if $\O$ is convex, then any Lipschitz constant in $\O_\d$ is also a Lipschitz constant in $\oli\O$.
\el
\demo
Set $u=u_{\O,h}$ and let $\d>0$ such that $u$ is $c$-Lipschitz in $\O_\d$ (i.e. Lipschitz with constant $c$). Let $z^1,z^2\in\oli\O$ with $\|z^2-z^1\|<\d/3$ and put $y=z^2-z^1$. Define the function $\Psi:\oli\O\to\R$ by
$$\Psi(z)=\left\{\begin{array}{ll} 
\max\{u(z),u(z+y)-c\|y\|\},& z\in\oli\O,\ dist(z,\bd\O)>\frac\d2\\ 
u(z),& z\in\oli\O,\ dist(z,\bd\O)\leq\frac\d2.\\
\end{array}\right.$$
Since $u$ is $c$-lipschitz in $\O_\d$, for every $z\in\O$ with $dist(z,\bd\O)=\frac\d2$
$$
|u(z+y)-u(z)|<c\|y\|,
$$
thus, by the Proposition \ref{lmax}, $\Psi$ is plurisubharmonic in $\O$ and therefore is a solution of the problem ($\star$). Using the maximality of $u$ we deduce that 
$$
u(z^2)-u(z^1)<c\|z^2-z^1\|.
$$
Interchanging $z^1$ and $z^2$ we obtain
$$
|u(z^2)-u(z^1)|<c\|z^2-z^1\|.
$$ 
Since every two points in $\oli\O$ can be connected by a path in $\oli\O$ made by segments of length less then or equal to $\d/3$, we get the conclusion. If $\O$ is convex, the path may be choosen a segment .
\enddemo
{\bf Proof of Proposition \ref{lip1}} Let $\tilde h:\oli\O\to\R$ be a $C^2$ extension of $h$ and $\phi$ a $C^2$ defining function for $\bd \O$ which is strongly plurisubharmonic. There exists a constant $a>0$ such that the function $\Phi^-=\tilde h+a\phi$ is plurisubharmonic and the function $\Phi^+=\tilde h-a\phi$ is plurisuperharmonic. Then we have $\Phi^-|_{\bd\O}=\Phi^+|_{\bd\O}=h$ and $\Phi^-\leq u\leq\Phi^+$. Since $\Phi^-$ and $\Phi^+$ are $C^2$, they are $c$-Lipschitz for some constant $k>0$.
Let $\d>0$, $z^1,z^2\in\oli\O$ such that $\|z^2-z^1\|=\d$ and put $y=z^2-z^1$. Let $z\in\O$ such that $dist (z,\bd\O)=\d$ and $\bar z\in\bd \O$ such that $\|z^0-z\|=\d$. Then, if $u=u_{\O,h}$, 
$$
u(z+y)\leq\Phi^+(z+y)\leq\Phi^+(z)+k\d\leq\Phi^+(\bar z)+2k\d\leq\Phi^-(z)+3k\d.
$$
In view of Proposition \ref{lmax}, the function
$$\Psi(z)=\left\{\begin{array}{ll} 
\max\{\Phi^-(z),u(z+y)-3k\d\},& z\in\oli\O,\ dist(z,\bd\O)>\d\\ 
\Phi^-(z),& z\in\oli\O,\ dist(z,\bd\O)\leq\d,\\
\end{array}\right.$$
is plurisubharmonic in $\O$, hence it is a solution of the problem ($\star$). Using the maximality of the Bremermann solution $u$, we get that 
$$u(z^2)-u(z^1)<k\|z^2-z^1\|$$
and interchanging $z^1$ and $z^2$ we obtain
$$
\vert u(z^2)-u(z^1)\vert<c\|z^2-z^1\|.
$$
Thus, $u_{\O,h}$ is Lipschitz on $\O_\d$, whence the conclusion in view of Lemma \ref{lip1}.
\enddemo

\section{Existence of continuous solutions for unbounded domains}\label{CPSS}
In this section we want to study the Bremermann-Dirichlet problem ($\star$), notably the existence of maximal solutions, when $\O$ is an unbounded domain in $\C^n$. It is worth to observe that the problem ($\star$) may not admit solutions even for very simple domains. Ideed, in \cite{ST99} it is proved that if $\Omega$ is the strongly convex paraboloid 
$$
\Omega=\left\{(z,w)\in\C^2: {\sf Im}\ w>|z|^2+({\sf Re}\ w)^2\right\},
$$
it is possible to construct a continuous function $h:\bd\Omega\rightarrow\R$ such that the only plurisubharmonic function $u$ on $\O$ which satisfy
$$
\limsup_{\z\rightarrow z}\ u(\z)\leq h(z),\ \forall z \in\Omega
$$
is the function $u\equiv-\infty$. The reason is that the function $h$ constructed there is negative in a "large part of $\bd\O$". Thus in the sequel we assume that the boundary value for the problem ($\star$) is non negative.

We start by proving the following existence theorem 
\bt \label{constr}
Let $\O\subset\C^n$ be a strongly pseudoconvex domain and $h:\bd\O\to [0,+\infty)$ a continuous function. Then the problem ($\star$) admits a solution $u$ which is continuous on $\oli\O$.
\et
\demo
For every $z\in \bd\O$ there exist a neighbourhood $U_z$ of $z$ and a biholomorphism $\phi_z: U_z\to U'_z$ such that $\phi_z(U_z\cap \bd\O)$ is strongly convex. We take two concentric balls $B'_z\subset C'_z\subset U'_z$ centered at $z'=\phi_z(z)$ and we define a continuous function $f'_z:U'_z\cap\phi(\bd\O)\to[0;1]$ such that
$$f'_z|_{\phi(\bd\O)\cap B'_z}\equiv 1,\ \ \ f'_z|_{U'_z\cap \phi(\bd\O)\ssmi C'_z}\equiv 0.$$
Put
$$S'_z=\{w\in U'_z\cap \phi(\bd\O):\ f'_z(w)>0\}.$$
If $S_z=\phi_{z_j}^{-1}(S'_z)$, then $\{S_z\}_{z\in \bd\O}$ is an open covering of $\bd\O$. From it, we extract a locally finite covering $\{S_{z_j}\}_{j\in J}$. If $f=f'\circ\phi_{z_j}$, we define the functions $\t_{z_j},\t:\bd\O\to [0,+\infty)$ by
$$\t_{z_j}(z)=\left\{\begin{array}{l}f(z),\ z\in S_{z_j}\\ 0,\ \ \ \ \ \rm{otherwise}\end{array}\right.$$
and $\t=\sum_{j\in J}\t_{z_j}$. Since $\{S_{z_j}\}_{j\in J}$ is locally finite, $\t<+\infty$, thus $\{\t_{z_j}/\t\}_{j\in J}$ is a partition of the unity. We define $h_{z_j}:\bd\O\to [0,+\infty)$ by $h_{z_j}=h\t_{z_j}/\t$, thus $\sum_{j\in J}h_{z_j}=h$. Putting $h'_{z_j}=h_{z_j}\circ \phi_{z_j}^{-1}$ and $G'_{z_j}=\phi_{z_j}(\oli\O)\cap C'_{z_j}$, we define $g'_{z_j}:\bd G'_{z_j}\to [0,+\infty)$ by
$$g'_{z_j}(z')=\left\{\begin{array}{l}h'_{z_j}(z'),\ z'\in S'_{z_j}\\ 0,\ \ \ \ \ \rm{otherwise}.\end{array}\right.$$
Notice that $G'_{z_j}$ may not be with $C^2$ boundary. For all that, we may assume $G'_{z_j}$ to be with $C^2$ boundary, using an approximation result proved in \cite{To83} (see also \cite{RS85}), which says that for every neighbourhood of the set of singularities of the boundary of a strongly convex analytic polihedron there exists a $C^2$ strongly convex subdomain whose boundary coincides with the boundary of the analytic polihedron outside of that neighbourhood. We consider the subset
$$D'_{z_j}=\{z'\in\oli G'_{z_j}:\ \exists L_{z'}\ni z' \text{ complex line s.t. }L_{z'}\cap \bd G'_{z_j}\subseteq\bd C'_{z_j}\}.$$
By the maximum principle, for every non-negative function $u'\in {\sf Psh}^{(s)}(G'_{z_j})$ with $u'|_{\bd G'_{z_j}}\leq g'_{z_j}$ we have $u'|_{D'_{z_j}}\equiv 0$. Thus, if $\G'_{z_j}$ is the solution for the Bremermann-Dirichlet problem for the domain $G'_{z_j}$ and the assigned values $g'_{z_j}$, then $\G'_{z_j}|_{D'_{z_j}}\equiv 0$, thus if $\G_{z_j}=\G'_{z_j}\circ\phi_{z_j}^{-1}$ and $D_{z_j}=\phi_{z_j}^{-1}(D'_{z_j})$, we have $\G_{z_j}|_{D_{z_j}}\equiv 0$. Since by construction $D'_{z_j}$ is a neighbourhood of $\bd C'_{z_j}\cap\phi_{z_j}(\oli\O)$ in $\oli G'_{z_j}$, $D_{z_j}$ is a neighbourhood of $C_{z_j}=\phi_{z_j}^{-1}(C'_{z_j})$ in $\oli G_{z_j}$, therefore the function $\psi_{z_j}:\oli\O\to\R$ given by
$$\psi_{z_j}=\left\{\begin{array}{l}\G_{z_j},\ \text{on }\oli\O\cap C_{z_j}\ssmi D_{z_j}\\ 0,\ \ \ \ \ \rm{otherwise}\end{array}\right.$$
is plurisubharmonic on $\O$. We define on $\oli\O$ the function $u$ by
$$u=\sum_{j\in J}\psi_{z_j},$$
which is a continuous solution of ($\star$).
\enddemo
\br
For every neighbourhood $U$ of $\bd\O$, $u$ may be constructed such that $u|_{\O\ssmi U}\equiv 0$.
\er
\bc
If $\O\subset\C^n$ is a strongly pseudoconvex domain, then there exists a function $\varrho\in C^0({\oli\O})\cap{\sf Psh}(\O)$ such that 
$\left\{\varrho=0\right\}=\bd\O$.
\ec
\demo
If we apply the previous theorem taking $h\equiv 1$, then the function $\phi:\oli\O\to\R$ given by $\phi=u-1$ satisfy the requirments. Indeed, assume that there exist $z^0\in\O$ such that $u(z^0)\geq1$. There exist $k\in\N$ and $J_1=\{j_1,\ldots,j_k\}\subset J$ such that $\psi_{z_j}(z^0)=0$ for every $j\in J\ssmi J_1$. Let $W_{z_j}=\rm{supp}\ \psi_{z_j}$ for $j\in J$. By construction, $\bd W_{z_j}=S_{z_j}\cup T_{z_j}$, such that $\psi_{z_j}|_{T_{z_j}}\equiv 0$ for all $j\in J_1$. On $\bigcup_{j\in J_1}W_{z_j}$ we define $u_1=\sum_{j\in J_1}\psi_{z_j}$. If we put $W=\rm{supp}\ u_1$, then $\bd W=S\cup T$, where $S=\bigcup_{j\in J_1}S_{z_j}$ and $T\subset\bigcup_{j\in J_1}T_{z_j}$ such that $\psi_{z_k}|_{T}\equiv 0$ for all $j\in J_1$. Therefore $u_1|_S\leq h\equiv1$ and $u_1|_T\equiv0$ which contradicts the maximum principle for the plurisubharmonic function $u_1$.
\enddemo
"Par abus de langage" also $\varphi$ will be called a {\it continuous defining function} for $\bd\O$.

\section{Maximal solutions for unbounded strongly convex domains}\label{msu}
The first case that we treat is the one of unbounded strongly convex domains of $\C^n$.

\subsection{Existence} Let $\O\subset\C^n$ be a strongly convex domain, $h:\bd\Omega\rightarrow[0,+\infty)$ a continuous function and consider the problem ($\star$) for the domain $\Omega$ and the function $h$. 

\begin{figure}[h]
\includegraphics{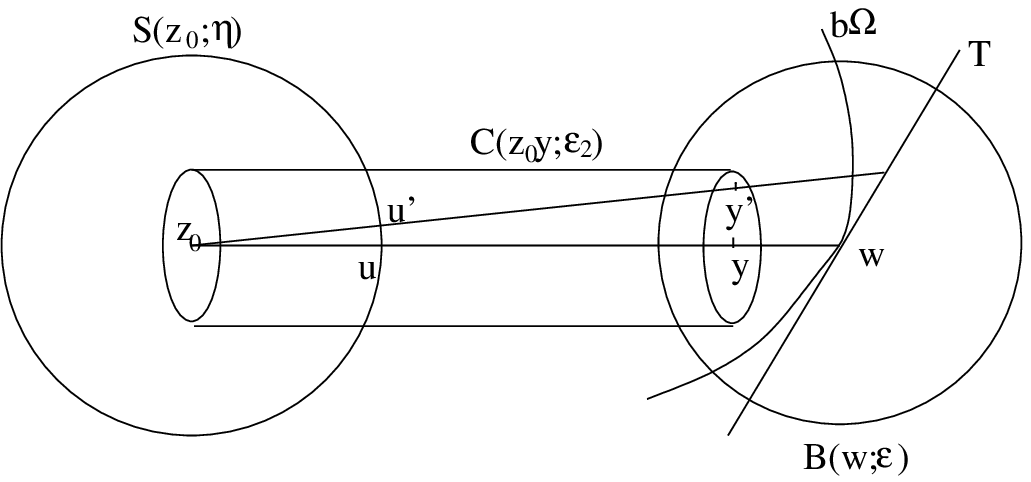}\caption{\Small }\label{fig4}
\end{figure}

We have the following intuitive fact
\bp \label{p2}
Let $\O$ be a convex unbounded domain in $\C^n$ with $C^2$ boundary. Then there exists $v\in\C^n\smi\{0\}$ such that for every $z\in\O$ and for every $k>0$, we have $z+kv\in\O$. Moreover, if $\O$ does not contain straight lines, there exists $v\in\C^n$ with the previous property and a unique complex hyperplane orthogonal to $v$ which is tangent to $\bd\O$.
\ep
\demo
Let $z^0\in\O$ and let $S(z^0;\eta)$ be a ball of center $z^0$. Suppose by contradiction that for every $z'\in S(z^0;\eta)$, the halfline starting in $z^0$ and passing by $z'$ intersects $\bd\O$ in some point $z$; we define $d(z')=\|z-z^0\|$. We will prove that $d:S(z^0;\eta)\to(0,+\infty)$ is a continuous function. Let $u\in S(z^0;\eta)$ and let $w$ be the intersection of the halfline starting in $z^0$ and passing by $u$ with $\bd\O$ and let $\e>0$. Let $H$ be the tangent hyperplane to $\bd\O$ in $w$. Let $y$ be a point on the segment of extremities $z^0$ and $w$, such that $\|y-w\|=\e_1$, with $\e_1<\e/2$ (see figure \ref{fig4}). There exists $\e_2<\e/2$ such that the cylinder $C(z^0y;\e_2)$ of radius $\e_2$ around the segment of extremities $z^0$ and $y$ to be contained in $\O$. Thus, the basis $B(y;\e_2)$ of the cylinder $C(z^0y;\e_2)$ which contains the point $y$ is included in the ball $B(w;\e)$. Decreasing, eventualy, $\e_2$, we may assume that for every point $y'\in B(y;\e_2)$, the halfline starting in $z^0$ and passing by $y'$ intersects $H$ inside the ball $B(w;\e)$. Taking the intersections of all such halflines with $S(z^0;\eta)$ we get a neighborhood $U$ of $u$ on $S(z^0;\eta)$ and taking the intersections of all such halflines with $\bd\O$, we get a neighborhood $W$ of $w$ in $\bd\O$ and $W\subset B(w;\e)$. Thus, for every $u'\in U$, $d(u')$ differs from $d(u)$ less than $\e$.

Since $\O$ is unbounded and convex, there exists a sequence $(u_\n)$ of points $u_\n\in S(z^0;\eta)$, such that $d(u_\n)\to+\infty$, which is a contradiction, since $d$ is continuous and $S(z^0;\eta)$ is compact. Therefore, there exists $v\in\C^n\smi\{0\}$ such that for every $k>0$, we have $z^0+kv\in\O$.

Let $z\in\O$. Suppose by contradiction that the halfline starting from $z$ and parallel with $v$ intersects $\bd\O$ in some point $w$. Since $\O$ is convex, the tangent hyperplane to $\bd\O$ in $w$ is not parallel with $v$, thus it will also intersect the halfline starting from $z^0$ and parallel with $v$, which contradicts again the convexity of $\O$.

Assume now that $\O$ contains no straight line and let
$$
C=\{u\in\C^n\smi\{0\}:\ \forall z^0\in\O,\ \forall k>0,\ z^0+ku\in\O\}.
$$
Since $\O$ is convex, if $u,v\in C$, then $\a u+\b v\in C$, for every $\a,\b\geq0$ and, since $\O$ contains no straight line, $C$ is included in an open halfspace. Let $v\in C$ such that the normal hyperplane $H$ to $v$ passing by the origin to be the boundaring hyperplane of an open halfspace including $C$. We may assume for simplicity that $v=(1,0,\ldots,0)$. Let
$$H_x=\{z=(x_1,y_1,\ldots,y_n)\in\C^n:\ x_1=x\}.$$
The first case we treat is when there exists $x\in\R$ such that $H_x\cap\O=\ES$. Let $l=\sup\{x\in\R:\ H_x\cap\O=\ES\}$.
Then $H_l\cap \bd\O\neq\ES$ and $H_l$ is tangent to $\bd\O$.

Now we treat together the subcase when $H_l\cap \bd\O=\ES$ (i.e.\ $H_l$ is asimpthotic to $\bd\O$) and the case when for every $x\in\R$, $H_x\cap\O\neq\Phi$. Let $w^0=(w_1,\ldots,w_n)\in\O$. We will prove that there exists $u\in\C^n\smi\{0\}$ orthogonal to $v$ such that for every $k>0$, we have $w^0+ku\in\O$. Assume by contradiction that the hyperplane $H_{w_1}$ passing by $w^0$ and orthogonal to $v$ intersects $\bd\O$ in a compact set $M$. For every $y\in M$, the tangent hyperplane $H_y$ to $\bd\O$ in $y$ intersects the line $\{w^0+kv:\ k\in\R\}$ in a point $p^y=(p_1,\ldots,p_n)$ with $p_1<w_1$. This implies that the set $\{z\in\O:\ z^1\leq w_1\}$ is included in the intersection of all halfspaces given by $H_y$ and $w^0$, thus is bounded. This is in contradiction with the assumption of the first case, as well as of the second case. Therefore there exists $u\in\C^n\smi\{0\}$ orthogonal to $v$ such that for every $k>0$, $w^0+ku\in\O$. We apply the previous proof to show that this property holds for every $z\in\O$ and thus $u\in C$. This is a contradiction, since $C$ is included in the open halfspace bounded by $H$.
\enddemo
As for the existence of a maximal solution the natural idea is to construct the maximal solution as limit of a sequence of Perron-Bremermann functions defined on bounded domains considering an appropriate exhaustive sequence of subdomains of $\O$. 

\begin{figure}[h]
\includegraphics{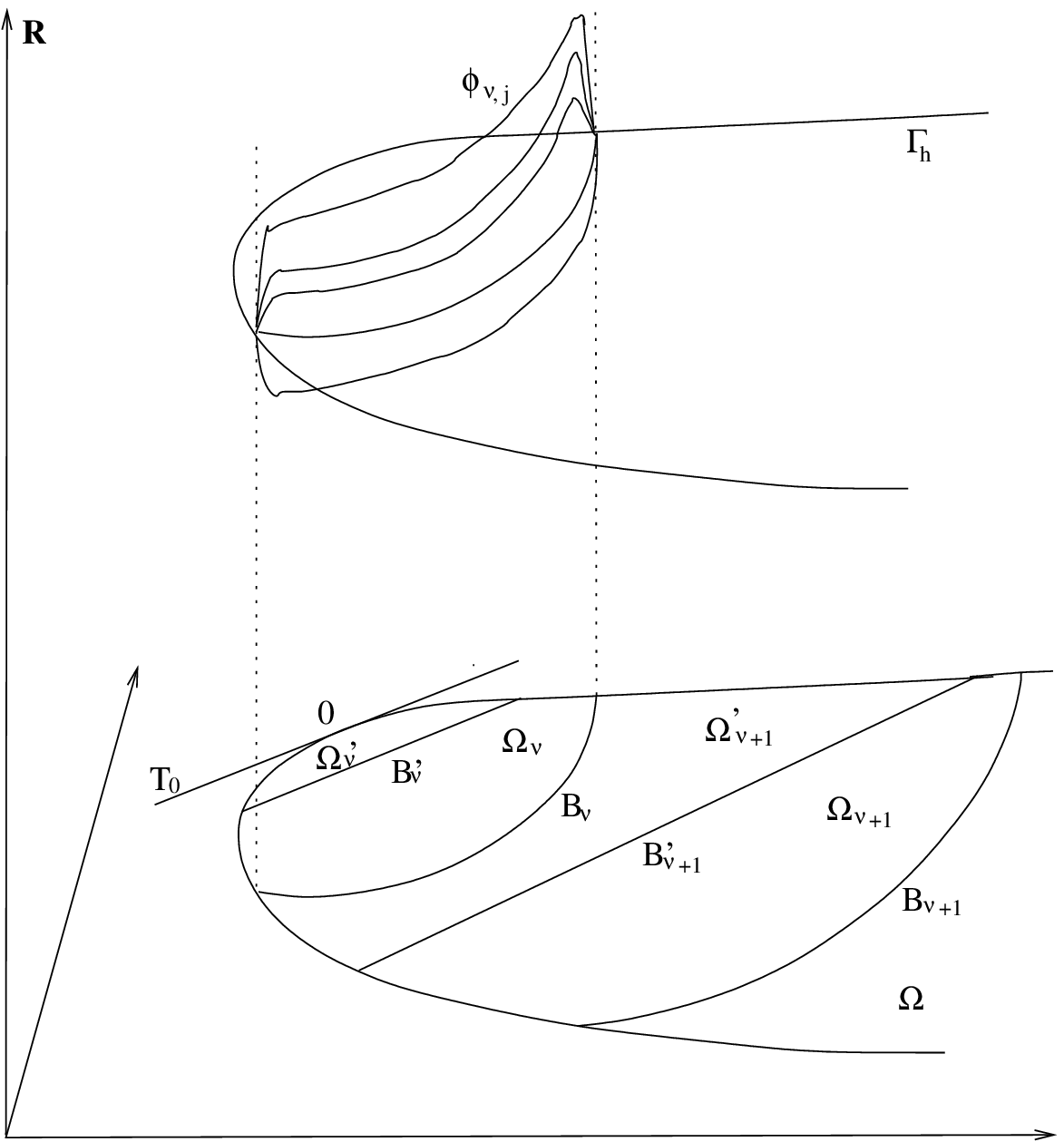}\caption{\Small }\label{fig1}
\end{figure}

We may assume that the vector given by the Proposition \ref{p2} to be $v=(1,0,\ldots,0)$ and let $k_0\in\R$ be such that $H_{0}=\{z=(x_1,y_1,\ldots,x_n,y_n)\in\C^n:\ x_1=k_0\}$ to be the unique hyperplane orthogonal to $v$ and tangent to $\bd\O$. We may assume that $k_0=0$. Let $(c_{\n})_{\n>0}$ and $(c'_{\n})_{\n>0}$ be increasing sequences in $\R$ such that 
$$\{z\in\O:\ x_1<c'_{\n}\}\subset B(0;c_\n)\subset \{z\in\O:\ x_1<c'_{\n+1}\},\ \forall n \in \N,$$
(see figure \ref{fig1}). Let 
$$\O'_{\n}=\{z\in\O:\ l(z)<c'_{\n}\}$$
and
$$\O_{\n}=\O\cap B(0;c_\n).$$
Notice that $\O_{\n}$ may not be with $C^2$ boundary. For all that, we may assume $\O_{\n}$ to be with $C^2$ boundary, using the same aproximation result as in the proof of the Theorem \ref{constr}. Let $A'=\bd\O'_{\n} \cap \bd\Omega$ and $B'_{\n}=\bd\O'_{\n} \cap H_{\n}$.
By Lemma \ref{l1}, if $u\in {\sf Psh}^{(s)}(\O_\n)$, then
$$\max_{z\in\oli\O'_\n}u(z)=\max_{z\in A'}u(z).$$
Let $A_{\n}=\bd\Omega_{\n} \cap \bd\Omega$ and $B_{\n}=\bd\Omega_{\n} \cap bB(0;c_\n)$. Let 
$$\{\phi:\bd\O_\n\to\R:\ \ \phi|_{A_\n}=h|_{A_\n}\}=\{\phi_{\n,j}\}_{j\in J}$$
be the set of all continuous extensions of $h|_{A_\n}$ to $\bd\O_\n$ (see figure \ref{fig1}). For every $\n \in \N$ and for every $j\in J$, let $\Phi_{\n,j}$ be the solution of the Bremermann-Dirichlet problem of the domain $\O_\n$ and function $\phi_{\n,j}$.
\bp\label{p88}
For every $\n\in\N$ and for every $z\in\O'_\n$ there exists $a_{z,\n}\in\R$, such that $\{\Phi_{\n,j}(z):\ j\in J\}=(-\infty,a_{z,\n})$ or $\{\Phi_{\n,j}(z):\ j\in J\}=(-\infty,a_{z,\n}]$.
\ep
\demo
For every $j\in J$ there exists a family of functions $\{\phi_{\n,t}\}_{t\in\R}$, $\phi_{\n,t}:\bd\O_\n\to\R$, for every $t\in\R$, which contains $\phi_{\n,j}$ as element, such that for every $z\in B_\n\smi A_\n$ we have:
\begin{itemize}
	\item [(i)] $\phi_{\n,t}(z)$ is continuous strongly increasing with respect to $t$;
	\item [(ii)] $\lim\limits_{t\to-\infty}\phi_{\n,t}(z)=-\infty$;
	\item [(iii)] $\lim\limits_{t\to+\infty}\phi_{\n,t}(z)=+\infty$.
\end{itemize}
By Proposition \ref{p1}, we get: if $t_1\leq t_2$, then $\phi_{\n,t_1}\leq\phi_{\n,t_2}$ and
$$\max_{z\in\oli\O_{\n,t}}|\Phi_{\n,t_1}-\Phi_{\n,t_2}|=\max_{z\in \bd\O_{\n,t}}|\phi_{\n,t_1}-\phi_{\n,t_2}|$$
for all $\n\in\N$ and $ t_1, t_2\in\R$. Since the function $\phi_{\n,t}(z)$ is continuous with respect to $t$, for every $z\in \bd\O_\n$, the function $\Phi_{\n,t}(z)$ is continuous with respect to $t$, for every $z\in \O_\n$, thus its image is an interval. Assume that $\{\Phi_{\n,t}(z):\ t\in\R\}$ is lower bounded in some $z^0\in\O'_\n$. Choosing $t_k=-k$, with $k\in\N$, we obtain a sequence of decreasing plurisubharmonic functions $(\Phi_{\n,k})_k$ and passing to the limit for $k\to\infty$, we obtain a plurisubharmonic function $\Psi_\n$ such that $\Psi_\n(z^0)>-\infty$. Let $z^1\in B_\n\smi A_\n$ and 
$\{\ l(z)=c\}$ be a complex hyperplane tangent in $z^1$ to $\O_\n$. Then there exists $\e_0>0$, such that $M=\{z\in\O_\n:\ l(z)=c+\e, |\e|<\e_0\}$ is a neighborhood of $z^1$ in $\O_\n$ which is union of domains in complex hyperplanes $M_\e=\{z\in\O_\n:\ l(z)=c+\e\}$ with boundary $bM_\e\subset (B_\n\smi A_\n)$
. Since $\Psi|_{M_\e}$ is plurisubharmonic and $\Psi|_{bM_\e}\equiv-\infty$, we get $\Psi|_M\equiv-\infty$, thus $\Psi\equiv-\infty$, contradiction.
Therefore, for $z\in\O_\n$ and for a family $\{\phi_{\n,t}\}_{t\in\R}$, the set $\{\Phi_{\n,t}(z):\ t\in\R\}$ is an interval which is not lower bounded. Thus the set $\{\Phi_{\n,j}(z):\ j\in J\}$ is union of intervals with the same property, which is also an interval not lower bounded.

Using Lemma \ref{l1}, we get that for every $z\in\O'_\n$, $\{\Phi_{\n,j}(z):\ j\in J\}$ is upper bounded by $\sup\limits_{z\in A'_\n} h(z)$.
\enddemo
\bp\label{p4}
Using the previous notations, if $\n_1<\n_2$, then for every $z\in\O'_{\n_1}$, $a_{z,\n_1}\geq a_{z,\n_2}$.
\ep
\demo
Assume by contradiction that there exists $z^0\in\O'_{\n_1}$ with $a_{z^0,\n_1}< a_{z^0,\n_2}$. Then there exists $j_2\in J$ such that $a_{z^0,\n_1}<\Phi_{\n_2,j_2}(z^0)$. Let $j_1\in J$, such that 
$$
\phi_{\n_1,j_1}(z)>\Phi_{\n_2,j_2}(z)
$$
for all $z\in B_{\n_1}$. We get the contradiction since $\Phi_{\n_1,j_1}$ is the solution of the Bremermann-Dirichlet problem for the function $\Phi_{\n_1,j_1}$, $\ \phi_{\n_1,j_1}(z)>\Phi_{\n_2,j_2}(z)$ on $\O_{\n_1}$, but $\Phi_{\n_1,j_1}(z^0)<\Phi_{\n_2,j_2}(z^0)$.
\enddemo
For every $\n\in\N$, let us define on $\oli\O'_\n$
$$\Phi_\n=\sup_{j\in J}\Phi_{\n,j}$$
and let
$$\Phi_\n^*(z)=\limsup_{\z\to z}\Phi_\n(\z)$$
be the upper semicontinuous envelope of $\Phi_\n$, which is a plurisubharmonic function.
\bp\label{p100}
Using the previous notations, for every $\n_0\in\N$, $\{\Phi_\n^*(z)\}_{\n\geq\n_0}$ is a decreasing sequence of functions, plurisubharmonic on $\O'_{\n_0}$ and for every $\n\in\N$, the function $\Phi_{\n+1}^*$ defined as before is continuous on $A_\n$ and satisfies $\Phi_{\n+1}^*|_{A_\n}=h|_{A_\n}$.
\ep
\demo
By the Proposition \ref{p4}, for every $\n_0\in\N$, the sequence $\{\Phi_\n(z)\}_{\n\geq\n_0}$ is decreasing on $\O'_{\n_0}$, thus passing to the $\limsup$ for $\z\to z$ in the points $z\in\O'_{\n_0}$, we get that $\{\Phi_\n^*(z)\}_{\n\geq\n_0}$ has the same property.

For every $\n\in\N$, $\Phi_\n$ was constructed as the $\sup$ of a family of continuous functions, therefore it is lower semicontinuous, thus for every $z\in A_\n$
\be\label{e1}
\Phi_{\n+1}^*(z)=\limsup_{\z\to z}\Phi_\n(\z)\geq\liminf_{\z\to z}\Phi_\n(\z)\geq\Phi_\n(z)=h(z)
\ee
Let $c:\bd\O'_{\n+1}\to\R$ be a continuous function such that
$$c|_{A_\n}=h|_{A_\n}=\Phi_{\n+1}|_{A_\n}$$
and
$$c|_{\bd\O'_{\n+1}\smi A_\n}\geq\Phi_{\n+1}|_{\bd\O'_{\n+1}\smi A_\n}.$$
There exists a continuous function $\G:\oli\O'_{\n+1}\to\R$, harmonic on $\O'_{\n+1}$, such that $\G|_{\bd\O'_{\n+1}}=c$. For every $j\in J$, $\Phi_{\n+1,j}$ is plurisubharmonic on $\O'_{\n+1}$ and thus, in particular, subharmonic. We have
$$\Phi_{\n+1,j}|_{\bd\O'_{\n+1}}\leq\Phi_{\n+1}|_{\bd\O'_{\n+1}}\leq c=\G|_{\bd\O'_{\n+1}},$$
for all $j\in J$, therefore 
$$
\Phi_{\n+1,j}|_{\O'_{\n+1}}\leq\G
$$
on $\O'_{\n+1}$ for all $j\in J$ . Thus 
$$
\Phi_{\n+1}|_{\O'_{\n+1}}\leq\G
$$
on $\O'_{\n+1}$. Passing to the $\limsup$ for $\z\to z$ in the points $z\in A_\n$ and using the continuity of $\G$, we get
\be\label{e2}
\Phi_{\n+1}^*|_{A_\n}\leq\G|_{A_\n}=h|_{A_\n}.
\ee
From (\ref{e1}) and (\ref{e2}) it follows $\Phi_{\n+1}^*|_{A_\n}=h|_{A_\n}$.
To prove that $\Phi_{\n+1}^*$ is continuous on $A_\n$, it suffices to show that it is lower semicontinuous on $A_\n$, since $\Phi_{\n+1}^*$ is upper semicontinuous. Indeed,
$$
\liminf_{\z\to z}\Phi_{\n+1}^*(\z)\geq\liminf_{\z\to z}\Phi_{\n+1}(\z)\geq\Phi_{\n+1}(z)=h(z)=\Phi_{\n+1}^*(z)
$$
for all $z\in A_\n$, since $\Phi_{\n+1}$ is lower semicontinuous.
\enddemo
For every $z\in\oli\O$ there exists $\n_0\in\N$, such that $z\in\oli\O_{\n_0}$. Let $\Phi:\oli\O\to\R$ be defined by
\be\label{Phi}
\Phi(z)=\inf_{\n\geq\n_0}\Phi_\n^*(z).
\ee
As a consequence of the way in which it was constructed, the function $\Phi$ extends plurisubharmonicaly the boundary values prescribing function $h$ and it is maximal among the functions with the same property. So we can conclude the section with the following
\bt\label{p6}
If $\O\subset\C^n$ is an unbounded strongly convex domain and $h:\bd\O\to [0,+\infty)$ a continuous function, then the Bremermann-Dirichlet problem ($\star$) admits a (non negative) maximal solution $\Phi$, i.e.\ if $\Phi'$ is another solution of ($\star$), then $\Phi'\leq\Phi$.
\et
\demo
We consider the previous construction of the function $\Phi$ from (\ref{Phi}). Since $h$ is non negative, $\Phi_\n$ is non negative for every $\n\in\N$, thus $\Phi$ is non negative. $\Phi$ was constructed as the infimum of a family of upper semicontinuous functions, thus it is upper semicontinuous. Since for every $\n\in\N$, $\Phi_{\n+1}^*|_{A_\n}=h|_{A_\n}$, it is clear that $\Phi|_{\bd\O}=h$.

Let $\Phi'\in {\sf Psh}^{(s)}(\O)$ be a function with $\Phi'|_{\bd\O}\leq h$, let $z\in\oli\O$ and let $\n\in\N$ such that $z\in\oli\O'_\n$. Let $j\in J$ such that $\Phi'|_{B_\n}\leq \phi_{\n,j}|_{B_\n}$. Then $\Phi'|_{\O'_\n}\leq \Phi_{\n,j}|_{\O'_\n}$, thus $\Phi'|_{\O'_\n}\leq \Phi_{\n}|_{\O'_\n}$ and $\Phi'|_{\O'_\n}\leq \Phi_{\n}^*|_{\O'_\n}$. Since the inequality holds for every $\n\in\N$ with $z\in\oli\O'_\n$, we get $\Phi'|_{\O'_\n}\leq \Phi|_{\O'_\n}$.
\enddemo
\br
The function $\Phi$ from the previous theorem is locally maximal.
\er
Indeed, if $G\subset\O$ in an open subset and $v\in {\sf Psh}^{(s)}(G)$ such that $v|_{\bd G}\leq \Phi|_{\bd G}$, we may consider the function $\Phi':\oli\O\to\R$ given by
$$\Phi'(z)=\left\{\begin{array}{ll} \max\{\Phi(z),v(z)\},& z\in\oli G\\ \Phi(z),& z\in\oli\O\ssmi G\end{array}\right.$$
By Proposition\ref{lmax}, $\Phi'$ is plurisubharmonic, therefore it is a solution of ($\star$). Using the maximality of $\Phi$, $\Phi'\leq\Phi$, therefore $v\leq \Phi$ on $G$.
\subsection{Continuity }\label{CMPS}
\bp\label{6789}
Let $\O\subset\C^n$ be an unbounded strongly convex domain and $h:\bd\O\to[0,+\infty)$ a continuous function. Then the maximal solution of the Bremermann-Dirichlet problem ($\star$) is continuous on $\bd\O$.
\ep
\demo
We consider the previous construction of the function $\Phi$ from (\ref{Phi}) and we have to prove that it is continuous on $\bd\O$. Since $\Phi$ is upper semicontinuous, it remains to prove its lower semicontinuity. Let $z^0\in \bd\O$.
By Theorem \ref{constr}, there exists a solution $u$ of the problem ($\star$) which is continuous on $\oli\O$. Let $\Psi:\oli\O\to\R$ be given by $\Psi=\max\{\Phi,u\}$. Then $\Psi$ is a solution of the problem ($\star$), thus by the Proposition \ref{p6} it follows that $\Phi\geq\Psi$ on $\oli\O$, which implies
$$\liminf_{z\to z^0}\Phi(z)\geq\liminf_{z\to z^0}\Psi(z)\geq\liminf_{z\to z^0}\Psi_1(z)=\Psi_1(z^0)=h(z^0)=\Phi(z^0).$$
\enddemo
\br
Continuity of  $h$ can be relaxed. Indeed, if $h$ is only upper semicontinuous, we can approximate it from above by a decreasing sequence of continuous functions $h_\n$. Thus we obtain a decreasing sequence of plurisubharmonic functions $\Phi_\n$ converging to a limit $\Phi$, which is plurisubharmonic and for each $z\in\bd\O$ such that $\limsup_{z\to z^0} h(z)=h(z^0)$ we have $\limsup_{z\to z^0} \Phi(z)=h(z^0)$.
\er
\bp
If $\O\subset\C^n$ is an unbounded strongly convex domain, $h:\bd\O\to [0,+\infty)$ a continuous function, $D\sbs\O$ is a strongly pseudoconvex subdomain such that the maximal solution $\Phi$ of the Bremermann-Dirichlet problem ($\star$) is continuous on $\bd D$, then $\Phi$ is continuous on $\oli D$.
\ep
\demo
Let $u$ denote the solution of the Bremermann-Dirichlet problem for the domain $D$ and given boundary $\Phi|_{\bd D}$. It is clear that $\Phi|_{D}\leq u$. Define 
$$\Psi(z)=\left\{\begin{array}{ll} u(z),& z\in\oli D\\ \Phi(z),& z\in\oli\O\ssmi\oli D,\end{array}\right.$$
By Proposition \ref{lmax}, the function $\Phi$ is plurisubharmonic, thus it is a solution of the problem ($\star$), therefore, by Theorem \ref{p6} we get $\Phi\geq\Psi$ and in particular $\Phi|_{\oli D}\geq\Psi|_{\oli D}=u$. Thus, $\Phi|_{\oli D}=u$, which is a continuous function.
\enddemo
The continuity of the maximal plurisubharmonic solution $\Phi$ on $\oli\O$ is more involved. It heavily depends on the behaviour of the boundary value $h$. In the rest of this section we treat some particular case which, in particular, allows us to prove continuity for bounded boundary values.

Let 
$$
M=\{x_1:\ z=(x_1,y_1,\ldots,y_n)\in \bd\O\}
$$
and $g:M\to\R$ be defined by
\be\label{e4}
g(x)=\max\{h(z):\ x_1\leq x\};
\ee
then
\be\label{e41}
g(x)=\max\{\Phi(z):\ x_1\leq x\}
\ee
by Lemma \ref{l1}.
 
 We have
\bl\label{p9}
If
$$
\lim_{x\to+\infty}\frac{g(x)}x=0,
$$
then the function $\Phi$ defined in (\ref{Phi}) is continuous.
\el
\demo
Let $z^0\in\O$, $z^0=(\xi_1,\z_1,\ldots,\xi_n,\z_n)$ and let $\e>0$. We define $R:\C^n\to\R$ by
$$R(z)=-\frac\e{2\xi_1}x_1-\frac\e2.$$
Since
$$\lim_{x\to+\infty}\frac{g(x)}x-\frac\e{2\xi_1}=-\frac\e{2\xi_1}<0,$$
we get
$$\lim_{x\to+\infty}g(x)-\frac\e{2\xi_1}x=-\infty,$$
thus there exists $\n\in\N$ such that 
\be\label{e3}
\Phi(z)+R(z)<0,\ \forall z=(x_1,y_1,\ldots,y_n)\in\O\smi\O'_{\n-2}.
\ee
Since $\Phi$ is continuous on $A'_\n$, for every $z\in A'_\n$ there exists $\d_z>0$ such that for every $z'\in B(z;2\d_z)$,
$$|\Phi(z)-\Phi(z')|<\frac\e4.$$
$\bigcup_{x\in A'_\n}B(z;\d_z)$ is an open covering of the compact set $A'_\n$, thus there exists a finite subcovering $\bigcup_{i\in I}B(z_i;\d_{z_i})$. Let $\d=\min\{\d_{z_i}: i\in I\}$. We claim that
$$
 \Phi(z)\geq\Phi(z^0)-\e,
$$
for all $z\in B(z^0;\d/3)$, i.e.\ $\Phi$ is lower semicontinuous in $z^0$, thus continuous. Let $z'\in B(z^0;\d/3)$ and $y=z^0-z'$.
Let $\Psi:\oli\O\to\R$ be given by
$$\Psi(z)=\left\{\begin{array}{ll} 
\max\{\Phi(z),\Phi(z+y)+R(z)\},& z\in\oli\O'_{\n-1},\ {\rm dist}(z,\bd\O)>\frac\d3\\ 
\Phi(z),& z\in\oli\O'_{\n-1},\ {\rm dist}(z,\bd\O)\leq\frac\d3\\
\Phi(z),& z\in\oli\O\smi\oli\O'_{\n-1}\\
\end{array}\right.$$
We will show that $\Psi$ is plurisubharmonic on $\O$. Since $\Phi$ is plurisubharmonic on $\O$, it remains to prove the plurisubharmonicity of $\Psi$ in a neighborhood of points $z$ with ${\rm dist}(z,\bd\O)=\d/3$ or $z\in B'_{\n-1}$. Let $z\in\O$ with ${\rm dist}(z,\bd\O)\in(\d/3,2\d/3)$. Then there exists $z''\in \bd\O$ with $\|z-z''\|<2\d/3$ and thus $\|(z+y)-z''\|<\d$. There exists $i\in I$ such that $\|z''-z_i\|<\d_{z_i}$, therefore $z,z+y\in B(z_i;2\d_{z_i})$, thus $$|\Phi(z)-\Phi(z_i)|<\e/4,\ \ \ |\Phi(z+y)-\Phi(z_i)|<\e/4,$$
thus
$$\Phi(z)>\Phi(z+y)-\e/2.$$
Therefore, in a neighborhood of points $z$ with ${\rm dist}(z,\bd\O)=\d/3$, $\Psi$ is equal to $\Phi$. Now, let $z\in B'_{\n-1}$ with ${\rm dist}(z,\bd\O)>\d/3$. Then, decreasing eventualy $\d$, $z+y\in\O'_\n\smi\O'_{\n-2}$ and using the relation (\ref{e3}) and the uniformly continuity of the function $R$ on $\oli\O'_\n\smi\O'_{\n-2}$, we get that
$$\Phi(z+y)+R(z)<0\leq\Phi(z)$$
in a neighborhood of $z$. Therefore, in a neighborhood of $z$, $\Psi$ is equal to $\Phi$.

By the Proposition \ref{p6}, we get that $\Phi\geq\Psi$. Taking $z=z'$, we get that $\Phi(z')\geq\Psi(z')\geq\Phi(z'+y)-\e=\Phi(z^0)-\e$.
\enddemo
As a immediate consequence, we get the following
\bp\label{2345}
Let $\O\subset\C^n$ be an unbounded strongly convex domain and $h:\bd\O\to\R$ a bounded continuous function. Then the maximal solution of the Bremermann-Dirichlet problem ($\star$) is continuous on $\oli\O$.
\ep
\br
Under the assumpsion of the previous proposition, (as well as for other less restrictive hypothesis which guarantee the continuity of the maximal solution, see Subsection \ref{se}), the maximal solution of the Bremermann-Dirichlet problem ($\star$) satisfies the distributional definition for the homogeneous complex Monge-Amp\`ere equation introduced by Bedford and Taylor (see \cite{BT76}).
\er
Some variations on the theme are contained in the next subsection.
\subsection{Some examples}\label{se}
\bl\label{p11}
Let $\O$ be a strongly convex domain in coordinates choosen like before. Let $m\in\N$, let
$$K=\{z=(x_1,y_1,\ldots,y_n)\in\C^n:\ m(2m-1)y_1^2-x_1^2<0,\ x_1>0\}$$
and assume that there exists $u=(\a,\b,0,\ldots,0)\in\C^n$ such that $\O\subseteq K-u$.  Using the previous notations, if
$$\lim_{x\to+\infty}\frac{g(x)}{x^{2m-1}}=0,$$
then $\Phi$ is continuous.
\el
\demo
The proof is analogue with the proof of the Lemma \ref{p9}. Let $a>0$. We put $$k=\frac\e{2(\xi_1+\a+a)^{2m}-2m(2m-1)(\z^1+\b)^2(\xi_1+\a+a)^{2m-2}}$$
and we define $R:\C^n\to\R$ by
$$R(z)=k(m(2m-1)(y_1+\b)^2(x_1+\a+a)^{2m-2}-(x_1+\a+a)^{2m})-\frac\e2.$$
The function $R$ has the following properties:
\begin{itemize}
	\item [(i)] $R\in {\sf Psh}(\C^n)$;
	\item [(ii)] $R(z)<-\e/2,\ \forall z\in K-u$;
	\item [(iii)] $R(z^0)=-\e$.
	\item [(iv)] $R(z)<-k(2ax_1^{2m-1}-a^2x_1^{2m-2}),\ \forall z=(x_1,y_1,\ldots,y_n)\in K-u$;
\end{itemize}
Since
$$\lim_{x\to+\infty}\frac{g(x)}{x^{2m-1}}-\frac{k(2ax^{2m-1}-a^2x^{2m-2})}{x^{2m-1}}=-2ak<0,$$
we get
$$\lim_{x\to+\infty}g(x)-k(2ax^{2m-1}-a^2x^{2m-2})=-\infty,$$
therefore, using the property (iv) and the definition of $g$ by relation (\ref{e4}), there exists $\n\in\N$ such that 
\be\label{e5}
\Phi(z)+R(z)<0,
\ee
for all $z\in\O\smi\O'_{\n-2}$.The rest of the proof goes as in Lemma \ref{p9}.
\enddemo
\bp
Let $\O\subset\C^n$ be an unbounded strongly convex domain in coordinates choosen like before and $h:\bd\O\to [0,+\infty)$ a continuous function. If there exist $a,b,c\in\R$ such that $\O\subseteq P$, where
$$P=\{z\in\C^n:\ ay_1^2+by_1+c<x_1\}$$
and there exists $m_0\in\N$ such that
$$\lim_{x\to+\infty}\frac{g(x)}{x^{m_0}}=0,$$
where $g$ is given by (\ref{e4}), then the maximal solution $\Phi$ of the Bremermann-Dirichlet problem ($\star$) is continuous.
\ep
\demo
Let $m\in\N$, such that $2m-1>m_0$ and let
$$K=\{z=(x_1,y_1,\ldots,y_n)\in\C^n:\ m(2m-1)y_1^2-x_1^2<0,\ x_1>0\}.$$
There exists $u=(\a,\b,0,\ldots,0)\in\C^n$ such that $P\subset K-u$. By the Lemma \ref{p11} it follows that $\Phi$ is continuous.
\enddemo
\bp\label{p12}
Let $\O\subset\C^n$ be an unbounded strongly convex domain in coordinates choosen like before and $h:\bd\O\to [0,+\infty)$ a continuous function. Let $a>0$ and
$$K=\{z\in\C^n:\ y_1\in(-\sqrt{2}/a,\sqrt{2}/a)\}.$$
If $\O\subseteq K$ and there exists $\a<a$ such that
$$\lim_{x\to+\infty}\frac{g(x)}{e^{\a x}}=0,$$
where $g$ is given by (\ref{e4}), then the maximal solution $\Phi$ of the Bremermann-Dirichlet problem ($\star$) is continuous.
\ep
\demo
We take
$$k=\frac{\e e^{-\a\xi_1}}{2\a^2\z_1^2-4}$$
and
$$R(z)=ke^{\a x}(\a^2y^2-2)-\e/2.$$
The condition (iv) becomes
$$R(z)<2k\left(\frac{\a^2}{a^2}-1\right)e^{\a x_1}$$
for all $z=(x_1,y_1,\ldots,y_n)\in K$ and since
$$\lim_{x\to+\infty}\frac{g(x)}{e^{\a x}}+2k\left(\frac{\a^2}{a^2}-1\right)=2k\left(\frac{\a^2}{a^2}-1\right)<0,$$
we get
$$\lim_{x\to+\infty}{g(x)}+2ke^{\a x}\left(\frac{\a^2}{a^2}-1\right)=-\infty$$
From now on the proof is similar to that of Lemma \ref{p11}.
\enddemo
\section{The case of unbounded strongly pseudoconvex domains}\label{9999}
The construction and the arguments used in the case of strongly convex domains in order to obtain the maximal plurisubharmonic function and to get its continuity on $\bd\O$, remain true in the case of the strongly pseudoconvex domains except for the "preparation theorem" and the construction of the "exhaustive" sequence of hyperplanes which have compact intersection with $\bd\O$. The domains $\O'_\n$, having the property that any plurisubharmonic function $u\in {\sf Psh}^{(s)}(\O_\n)$ riches its maximum on that part of the boundary of $\O'_{\n}$ which is containd in $\bd\O$ (denoted by $A'_{\n}$), were obtained in the strongly convex case as the intersection of $\O$ with a convenient sequence of parallel halfspaces. In the case of the strongly pseudoconvex domains it may be possible for any sequence of hyperplanes to contain elements which divide $\O$ in two subsets which are both not relatively compact. Therefore, the domains $\O'_\n$ will be constructed in a different way, namely as envelopes of plurisubharmonicity. 

In this section $\Omega\subset\C^n$ will denote a strongly pseudoconvex domain. Let $(\Sigma_\n)_{\n\in\N}$ be a sequence of compact subsets $\Sigma_\n\subset \bd\O$, such that $\Sigma_\n\subset\Sigma_{\n+1}$ for every $\n\in\N$ and
$$\bigcup_{\n\in\N}\Sigma_\n=\bd\O.$$
Passing, if necessary, to a subsequence, we may find a sequence of balls $B(0;c_\n)$ with increasing radius $c_\n$ such that if we denote as in the convex case
$$\O_{\n}=\O\cap B(0;c_\n),$$
$A_{\n}=\bd\Omega_{\n} \cap \bd\Omega$ and $B_{\n}=\bd\Omega_{\n} \cap bB(0;c_\n)$, then
$$A_{\n-1}\subset\Sigma_{\n}\subset A_{\n},$$
for every $\n\geq 1$. We define $\O'_\n$ as the plurisubharmonic envelope of $\Sigma_\n$
$$\O'_\n=\{z\in \oli\O_\n:\ f(z)\leq\max_{\z\in\Sigma_\n}f(\z),\ \forall f\in {\sf Psh}^{(s)}(\O_\n)\}$$

Differently from what happens in the case of strongly convex domains, a priori we may now have
\be\label{nestrict}
\bigcup_{\n\in\N}\Omega'_\n\subsetneq\Omega.
\ee
which would not guarantee $\Phi<+\infty$. A sufficient condition in order to have equality in (\ref{nestrict}) is that for every point $z\in\Omega$ there exists an analytic disc $D_z\ni z$ such that $D_z\cap \bd\Omega$ to be compact. A sufficient condition for having this was first pointed out by Lupacciolu \cite{Lu87}
in studying the extension problem for $CR$ functions in unbounded
domains. Precisely,
\begin{enumerate}
\item[($L$)] if $\oli \Omega^\infty$ denotes the closure of $\Omega \subset \C^n\subset \C\P^n$ in
$\C\P^n$, then there exists an algebraic hypersurface
$$Z=\{z\in\C^n:\ P(z)=0\}$$
such that $Z\cap\oli \Omega^\infty=\emptyset$.
\end{enumerate}
Equivalently
\begin{enumerate}
\item[($L'$)] there exists a polynomial $P\in \C[z_1,\ldots,z_n]$ such that
$$\Omega \subset\left\{z\in \C^n:\vert P(z)\vert^2>(1+\vert z\vert^2)^{{\rm deg} P}\right\}.$$
\end{enumerate}
\br
There exist unbounded strongly pseudoconvex domains $\O\subset\C^n$ which do not verify the condition ($L$).
\er
Indeed, if we consider a polynomial $Q\in \C[z_1,\ldots,z_n]$ and the domain
$$\O=\{z\in\C^n:\ \|z\|^2+\log \|Q(z)\|^2<1\},$$
then $\O$ is strongly pseudoconvex since its boundary is a level set of a strongly plurisubharmonic function. If $Z_Q$ is the zero set of $Q$ we have $Z_Q\subset\O$. Then for every polynomial $P\in \C[z_1,\ldots,z_n]$, the zero set $Z_P$ of $P$ and $Z_Q$ intersects in $\C\P^n$, i.e.\ condition ($L$) is not satsfied.
\br\label{DEF}
It is worthy observing that $\O$  has no defining function. Indeed,  assume that $\varrho$ is strongly plurisubharmonic in a neighbourhood of $\oli\O$ and such that $\O=\{\varrho<0\}$. Then, $\varrho|_{Z_Q}$ is bounded above, thus extendable on the projective closure $\widehat Z_Q$ of $Z_Q$ hence constant: a contradiction.
\er
\bt
If $\O\subset\C^n$ is an unbounded strongly pseudoconvex domain which verifies the condition ($L$) and $h:\bd\O\to [0,+\infty)$ a continuous function, then the Bremermann-Dirichlet problem ($\star$) admits a (non negative) maximal solution $\Phi$, i.e.\ if $\Phi'$ is an other solution of ($\star$), then $\Phi'\leq\Phi$. Moreover, the function $\Phi$ is continuous on $\bd\O$.
\et
\demo
The only thing we have to show in order to conclude the proof (by
using the methods of the previous section) is that, up to a
holomorphic change of coordinates and a holomorphic embedding $E:\C^n\rightarrow \C^N$, for every $z\in\Omega$, there exists a hyperplane passing by $z$ which intersects $\bd\Omega$ in a compact set.

Following \cite{Lu87}, we divide the proof in two steps.

\emph{Step 1}. $P$ linear. We consider
$\oli\O\subset\C\P^n=\C^n\cup\C\P^{n-1}_\infty$, which is disjoint
from $Z=\left\{P=0\right\}$. So we can consider new coordinates
of $\C\P^n$ in such a way that $Z$ is the $\C\P^{n-1}$ at
infinity. Now $\O$ is a relatively compact open set of
$(\C^n)'=\C\P^n\smi Z$ and $H_\infty=\C\P^{n-1}_\infty\cap(\C^n)'$ is a complex hyperplane
containing the topological boundary of $\bd\Omega$. For every $z\in\Omega$, the complex hyperplane passing by $z$ and parallel to $H_\infty$ intersects $\bd\Omega$ in a compact set $M_z$. Similarely to the case when $\O$ was strongly convex, we can define an exaustive sequence $\O'_\n$ of subsets of $\O$ such that every function $u\in {\sf Psh}^{(s)}(\O'_\n)$ reaches its maximum on $A'_\n=\bd\O'_\n\cap\bd\O$. For instance, we can take
$$\O'_\n=\bigcup_{\dist (z,H_\infty)>1/n}M_z.$$

\emph{Step 2}. $P$ generic. We use the Veronese map $V$ to embed
$\C\P^n$ in a suitable $\C\P^N$ in such a way that $V(Z)=L_0\cap
V(\C\P^n)$, where $L_0$ is a linear subspace. The Veronese map $V$
is defined as follows: let $d$ be the degree of $P$, and let
$$N\ =\ {{n+d}\choose{d}}-1.$$
Then $V$ is defined by
$$V(z)\ =\ V[z^0:\ldots:z_n]\ =\ [\ldots:w_I:\ldots]_{|I|=d},$$
where $w_I=z^I$. If $P=\sum_{|I|=d}\alpha_I z^I$, then $V(Z)=L_0\cap V(\C\P^n)$, where
$$L_0=\left\{\sum_{|I|=d} \alpha_I w_I = 0 \right\}.$$
Again we can change the coordinates such that $L_0$ is the $\C\P^{N-1}$ at infinity. We may now find the complex hyperplane as in Step 1.
\enddemo

\section{Existence of a pluriharmonic solution}
Similarly to the plurisubharmonic functions, one can define the plurisuperharmonic functions. Precisely, a plurisuperharmonic function $\Phi:\O\to\R\cup\{+\infty\}$ on a domain $\O\subseteq\C^n$ is a lower semicontinuous function such that its restriction to any complex line is superharmonic; equivalently, $\Phi$ is plurisuperharmonic if and only if $-\Phi$ is plurisubharmonic.

Obviously, we can obtain for plurisuperharmonic functions results which are similar to those we get in the previous sections for plurisubharmonic functions. Thus, if $\O\subset\C^n$ is an unbounded strongly convex domain and $h:\bd\O\to\R$ a continuous function, then the Bremermann-Dirichlet problem 
\begin{itemize}
	\item [($\star'$)] $\left\{\begin{array}{cl}(i)& u \text{ is lower semicontinuous in }\oli\O\\ (ii)&u\text{ is plurisuperharmonic in }\O\\ (iii)&u|_{\bd\O}=h\end{array}\right.$
\end{itemize}
(which for making sense, we assume $h$ upper bounded), admits a minimal solution $\chi$ which is continuous on $\bd\O$. Moreover, if $h$ is bounded, then $\chi$ is continuous on $\oli\O$.

For the solvability of both ($\star$) and ($\star'$) problems, we assume that $h$ is bounded. In this case, the respective solutions $\Phi$ and $\chi$ are continous on $\oli\O$, satisfy $\Phi\leq\chi$ and for every harmonic solution $\eta$ of the classical Dirichlet problem we have $\Phi\leq\eta\leq\chi$.

If there exists $z\in\O$ such that $\Phi(z)=\chi(z)$, then for the maximum principle applied to the plurisubharmonic function $\Phi-\chi$, we get $\Phi\equiv\chi$, thus there exists a pluriharmonic solution of the problem.

It is well known that boundary values of pluriharmonic functions must satisfy certain compatibility conditions. In \cite{DBT} it is constructed a differential operator of the third order $\oli\partial_{\rm b}\o$ such that: on a strongly pseudoconvex hypersurface, the traces $h$ in the sense of distributions of pluriharmonic functions are solutions of $\oli\partial_{\rm b}\o$; conversely, if $\Omega$ is simply connected and ${\rm b}\Omega$ is strongly pseudoconvex, every solution of the operator $\oli\partial_{\rm b}\o$ can be extended to a pluriharmonic function.

\section{A generalization for $q$-plurisubharmonic functions}
A $C^2$ function $u:U\to\R$ on an open set $U\subseteq\C^n$ is called $q$-plurisubharmonic, $0\leq q\leq n-1$, if its Levi form $L(u)$ has at least ($n-q$) non-negative eigenvalues. A generalization of the problem ($\star$) for $q$-plurisubharmonic functions was considered by Hunt and Murray in \cite{Hu78}. The authors gave a new definition for $q$-plurisubharmonic functions, which is also applicable to less regular functions. Precisely, a $q$-plurisubharmonic function $u:U\to\R\cup\{-\infty\}$ on an open set $U\subseteq\C^n$ is an upper semicontinuous function such that for every ($q+1$)-dimensional complex plane $L$ intersecting $U$, for every compact subset $K\subset U\cap L$ and for every superharmonic function $v$ on $\rm{int}\,K$, $C^2$ up to the boundary such that $v|_{\bd K}\geq u|_{\bd K}$, we have $v|_{K}\geq u|_{K}$. Observe that with this terminology, the $0$-plurisubharmonic functions are the plurisubharmonic functions. If $U\in\C^n$ is an open set, we denote by ${\sf Psh}^{(s)}_q(U)$ the set of upper semicontinuous functions in $\oli U$ and $q$-plurisubharmonic in $U$. 

Given a domain $\O\subset\C^n$ and a continuous function $h:\bd\O\to\R$, a generalized Dirichlet problem consists of finding a function $u:\oli\O\to\R$, such that
\begin{itemize}
	\item [($\star''$)] $\left\{\begin{array}{cl}(i)& u \text{ is upper semicontinuous in }\oli\O\\ (ii)&u\text{ is } q\text{-plurisubharmonic in }\O\\ (iii)&u|_{\bd\O}= h\end{array}\right.$
\end{itemize}

Let $\O\subseteq\C^n$ be a domain with $C^2$ boundary. Let $\varrho$ be a defining function for $\O$ in a neighbourhood $U$ of a point $p\in\bd \O$. If the Levi form $L_p(\varrho)$ of $\varrho$ at $p$ has at most $q$ negative eigenvalues, $\O$ is called strongly $q$-pseudoconvex at $p$. A bounded strongly $q$-pseudoconvex domain $\O$ admits a $q$-plurisubharmonic globaly defining function in a neighbourhood of $\oli\O$. Hunt and Murray proved that if $\O$ is a strongly $q$-pseudoconvex bounded domain and $2q<n$, then the problem ($\star''$) admits a maximal solution $u^q_{\O,h}$; this is the upper envelope of the functions $u$ which are upper semicontinuous on $\oli\O$, $q$-plurisubharmonic on $\O$ and such that $u|_{\bd\O}\leq h$. Moreover, $u^q_{\O,h}$ is continuous on $\oli\O$ and it is also $(n-q-1)$-plurisuperharmonic. The same is true in a slightly more general situation, namely for bounded strongly $r$-pseudoconvex domains and $q$-plurisubharmonic functions with $0\leq q\leq n-r-1$, as noticed by Slodkowski in \cite{Sl84}. In the same paper, the author call $q$-Bremermann a continuous function which is both $q$-plurisubharmonic and $(n-q-1)$-plurisuperharmonic and proved that for a continuous function $h$ and for each $0\leq q\leq n-1$ there exists at most one function $u$ continuous on $\oli\O$, $q$-Bremermann on $\O$ and such that $u|_{\bd\O}= h$. In particular, if $\O$ is a bounded strongly pseudoconvex domain and $h$ a continuous function, then for each $0\leq q\leq n-1$ the problem ($\star''$) has a unique $q$-Bremermann solution.

In this section we want to extend to unbounded domains the results concerning the $q$-plurisubharmonic functions and the problem ($\star''$), as we did in the previous sections for the plurisubharmonic functions and the problem ($\star$).
We assume again that the vector given by the Proposition \ref{p2} is $v=(1,0,\ldots,0)$ and we repeat for the $q$-plurisubharmonic functions the construction we made in Section \ref{msu} for plurisubharmonic functions. The Propositions \ref{lmax}, \ref{p1}, \ref{p88}, \ref{p4} hold for $q$-plurisubharmonic functions with similar proofs. Lemma \ref{l1}, which was used to prove the upper boundedness of a family of plurisubharmonic functions, holds also for $q$-plurisubharmonic functions, since the harmonic function $\phi\equiv M$ is in particular superharmonic. The Proposition \ref{p100} holds for $q$-plurisubharmonic functions with a slightly different proof, namely the function $\Gamma$ must be choosen plurisuperharmonic instead of harmonic, so we will take $\Gamma$ to be the solution of the Bremermann-Dirichlet problem for plurisuperharmonic functions for the domain $\O'_{\nu +1}$ and boundary value $c$. All this results allow us to state the following
\bt\label{qps}
If $\O\subset\C^n$ is an unbounded strongly convex domain and $h:\bd\O\to [0,+\infty)$ a continuous function, then for every $1\leq q\leq n-1$ the generalized Dirichlet problem ($\star''$) admits a (non negative) maximal solution $\Phi$, i.e.\ if $\Phi'$ is an other solution of ($\star''$), then $\Phi'\leq\Phi$. Moreover, $\Phi$ is continuous on $\bd\O$.
\et
\demo
The first part of the proof follows as for plurisubharmonic functions. To prove the continuity of the function $\Phi$ on $\bd\O$, we observe that the function furnished by the Theorem \ref{constr} is in particular $q$-plurisubharmonic.
\enddemo
\bt\label{conq}
Let $\O\subset\C^n$ be an unbounded strongly convex domain in coordinates choosen like before and $h:\bd\O\to [0,+\infty)$ a continuous function. If one of the following conditions hold:
\begin{itemize}
	\item [(i)] $\lim\limits_{x\to+\infty}\frac{g(x)}x=0$ or in particular $h$ is bounded;
	\item [(ii)] there exist $a,b,c\in\R$ such that $\O\subseteq P$, where
$$P=\{z\in\C^n:\ ay_1^2+by_1+c<x_1\}$$
and there exists $m_0\in\N$ such that $\lim\limits_{x\to+\infty}\frac{g(x)}{x^{m_0}}=0$;
	\item [(iii)] there exist $a>0$ such that $\O\subseteq K$, where
$$K=\{z\in\C^n:\ y_1\in(-\sqrt{2}/a,\sqrt{2}/a)\}$$
and there exists $\a<a$ such that $\lim\limits_{x\to+\infty}\frac{g(x)}{e^{\a x}}=0$,
\end{itemize}
where $g$ is given by (\ref{e4}), then the maximal solution $\Phi$ of the generalized Dirichlet problem ($\star''$) is continuous in $\oli\O$ and ($n-q-1$)-plurisuperharmonic in $\O$.
\et
\demo
The proof of the continuity is similar with the case of plurisubharmonic functions. In this situation we use a result of Sadullaev (see \cite{Sl84}) which says that the sum between a $q$-plurisubharmonic function and a $r$-plurisubharmonic function is ($q+r$)-plurisubharmonic and in the particular case $r=0$ the sum between a $q$-plurisubharmonic function and a plurisubharmonic function is $q$-plurisubharmonic.

The proof of the fact that the maximal solution $\Phi$ is ($n-q-1$)-plurisuperharmonic in $\O$ is similar to the case of bounded domains. Let $L$ be a ($q+1$)-dimensional complex plane intersecting $\O$, $K\subset \O\cap L$ a compact subset and $v$ a subharmonic function on $\rm{int}\,K$, $C^2$ up to the boundary, such that $v|_{\bd K}< \Phi|_{\bd K}$. We may assume that $K$ is a ball and that the complex plane $L$ is given by $z_1=0$, $z_2=0,\ldots,z_q=0$, thus $v(z)$ is actually $v(z_{q+1},\ldots,z_n)$. Let $\widetilde v(z_1,z_2,\ldots,z_n)$ be given by
$$\widetilde v(z_1,z_2,\ldots,z_n)=v(z_{q+1},\ldots,z_n)-k(z_1\bar z_1+\ldots +z_q\bar z_q),$$
where $k>0$ is choosen sufficiently large such that $\widetilde v<\Phi$ on $\bd U$, where $U\subset\O$ is still to be choosen strongly pseudoconvex domain with $C^2$ boundary and with closure in $\O$ whose intersection with $L$ is $K$. Since the Levi form of $\tilde v$ has ($n-q$) non-negative eigenvalues, $\tilde v$ is $q$-plurisubharmonic on $U$ if for instance $U$ is choosen an ellipsoid which has $L$ as symmetry plane. Since $\widetilde v<\Phi$ on $\bd U$, using the result similar to the Proposition \ref{lmax} for $q$-plurisubharmonic functions we get that the function $\Psi:\oli \O\to\R$ given by
$$\Psi(z)=\left\{\begin{array}{ll} \max\{\Phi,\widetilde v\},& \text{on } U\\ \Phi,& \text{on } \O\ssmi\oli U,\end{array}\right.$$
is a solution of the problem ($\star''$). Using the maximality of the function $\Phi$, we get that $\tilde v|_{U}< \Phi|_{U}$ and thus $v|_{K}=\tilde v|_{K}< \Phi|_{K}$.
\enddemo

\end{document}